\documentclass[12pt]{article}
\usepackage{latexsym,color,amsmath,amsthm,amssymb,amscd,amsfonts}

\newtheorem{lemma}{Lemma}[section]
\newtheorem{proposition}[lemma]{Proposition}
\newtheorem{remark}[lemma]{Remark}
\newtheorem{theorem}[lemma]{Theorem}
\newtheorem{definition}[lemma]{Definition}
\newtheorem{corollary}[lemma]{Corollary}

\newtheorem*{remark*}{Remark}

\def\R{{\mathbb R}}
\def\sph{S^{2n-1}}
\def\eps{\varepsilon}

\makeatletter \@addtoreset {equation}{section}

\renewcommand\theequation
  {\ifnum \c@subsection>\z@ \arabic{section}.\arabic{subsection}.\arabic{equation}
  \else \arabic{section}.\arabic{equation} \fi}
\makeatother

\begin{document}
\title {A Brunn-Minkowski Inequality for Symplectic Capacities of Convex
Domains}
\author{Shiri Artstein-Avidan, \ Yaron Ostrover\thanks{The first and
second named authors were both partially supported by BSF grant
no. 2006079. \hfill 
This first named author was partially supported by the ISF grant No.
865/07, and the second named author was partially supported by NSF
grant DMS-0706976.}} \maketitle \noindent {\it {\bf Abstract:} In
this work we prove a Brunn-Minkowski-type inequality in the context
of symplectic geometry and discuss some of its applications. }

\section{Introduction and Results}

In this note we examine the classical Brunn-Minkowski inequality in
the context of symplectic geometry. Instead of considering volume,
as in the original inequality, the quantity we are interested in is
a symplectic capacity, given by the minimal symplectic area of a
closed characteristic on the boundary of a convex domain. To explain
the setting, the main results, and their significance, we start with
an introduction.

\subsection{The Brunn-Minkowski inequality} \label{BM-section}

Denote by ${\cal K}^n$ the class of convex bodies in ${\mathbb
R}^n$, that is, compact convex sets with non-empty interior. The
Brunn-Minkowski inequality, in its classical formulation, states
that if $K$ and $T$ are in ${\cal K}^n$, then
\[  ({\rm Vol}( K + T)
)^{\frac 1 n} \geq ({\rm Vol}(K))^{\frac 1 n} + ({\rm
Vol}(T))^{\frac 1 n},\] where ${\rm Vol}$ denotes the
$n$-dimensional volume (i.e, the Lebesgue measure) and the
Minkowski sum of two bodies is defined by
$$ K + T = \{  x + y \ : \ x \in K, \ y\in T \}.$$
Moreover, equality holds if and only if $K$ and $T$ are homothetic,
or in other words, coincide up to translation and dilation.

The Brunn-Minkowski inequality is a fundamental result in convex
geometry and has innumerable  applications, the most famous of which
is probably a simple proof of the isoperimetric inequality. We
recall that in fact it is known that the Brunn-Minkowski inequality
holds for any two measurable sets (while the equality condition
requires convexity, or some special form of non-degeneracy). The
inequality is connected with many other important inequalities such
as the isoperimetric inequality, the Sobolev and the Log-Sobolev
inequalities, and the Pr\'{e}kopa-Leindler inequality. Moreover, the
Brunn-Minkowski inequality has diversified applications in analysis,
geometry, probability theory, information theory, combinatorics,
physics and more. We refer the reader to \cite{Gar} for a detailed
survey on this topic.

The importance of the Brunn-Minkowski inequality has led to efforts
of finding analogous inequalities in other areas of mathematics, and
recently, inequalities of Brunn-Minkowski type were proved for
various well known functionals other than volume. Two examples of
these functionals, which are related with calculus of variation and
with elliptic partial differential equations, are the first
eigenvalue of the Laplace operator~\cite{BL} and the electrostatic
capacity~\cite{B}. In this note we concentrate on a symplectic
analogue of the inequality. To explain it, we turn now to the
framework of symplectic geometry.

\subsection{Symplectic Capacities} \label{sect-sympac}

Consider the $2n$-dimensional Euclidean space ${\mathbb R}^{2n}$
with the standard linear coordinates $(x_1,y_1, \ldots,x_n,y_n)$.
One equips this space with the standard symplectic structure
$\omega_{st} = \sum_{j=1}^n dx_j \wedge dy_j$, and with the standard
inner product $g_{st} = \langle \cdot,\cdot \rangle$.
Note that under the identification 
between ${\mathbb R}^{2n}$ and $ {\mathbb C}^n$, these two
structures are the real and the imaginary parts of the standard
Hermitian inner product in ${\mathbb C}^n$, and $\omega(v,Jv) =
\langle v , v \rangle $, where $J$ is the standard complex structure
on ${\mathbb R}^{2n}$. Recall that a {\it symplectomorphism} of
${\mathbb R}^{2n}$ is a diffeomorphism which preserves the
symplectic structure i.e., $\psi \in {\rm Diff} ( {\mathbb R}^{2n}
)$ such that $\psi^* \omega_{st} = \omega_{st}$. In what follows we
denote by ${\rm Symp}({\mathbb R}^{2n})$ the group of all the
symplectomorphisms
 of ${\mathbb R}^{2n}$.

Symplectic capacities are symplectic invariants which, roughly
speaking, measure the symplectic size of subsets of ${\mathbb
R}^{2n}$. More precisely,

\begin{definition} \label{Def-sym-cap}
A symplectic capacity on $({\mathbb R}^{2n},\omega_{st})$ associates
to each  subset $U \subset {\mathbb R}^{2n}$ a number $c(U) \in
[0,\infty]$ such that the following three properties hold:
\begin{enumerate}
\item[(P1)] $c(U) \leq c(V)$ for $U \subseteq V$ (monotonicity)
\item[(P2)] $c \big (\psi(U) \big )= |\alpha| \, c(U)$ for  $\psi
\in {\rm Diff} ( {\mathbb R}^{2n} )$ such that $\psi^*\omega_{st} =
\alpha \, \omega_{st}$ (conformality) \item[(P3)] $c \big (B^{2n}(r)
\big ) = c \big (B^2(r) \times {\mathbb C}^{n-1} \big ) = \pi r^2$
(nontriviality and normalization),
\end{enumerate}
\end{definition}

\noindent where $B^{2k}(r)$ is the open $2k$-dimensional Euclidean
ball of radius $r$. Note that the third property disqualifies any
volume-related invariant, while the first two properties imply that
for two sets $U,V \subset {\mathbb R}^{2n}$, a necessary condition
for the existence of a symplectomorphism $\psi$ such that $\psi
(U_1) = U_2$ is that $c(U_1) = c(U_2)$ for each symplectic capacity
$c$.

A priori, it is unclear that symplectic capacities exist. The
first example of a symplectic capacity is due to Gromov~\cite{G}.
His celebrated non-squeezing theorem states that for $R
> r$ the ball $B^{2n}(R)$ does not admit a symplectic embedding
into the symplectic cylinder $Z^{2n}(r):= B^2(r) \times {\mathbb
C}^{n-1}$. This theorem led to the following definitions:

\begin{definition} The symplectic radius of a non-empty set
$U \subset {\mathbb R}^{2n}$ is
$$ c_B(U) := \sup \left \{\pi r^2 \, | \,
\ There \ exists \  \psi \in {\rm Symp}({\mathbb R}^{2n}) \ with \
\psi \left (B^{2n}(r) \right ) \subset U \right  \}.$$ The
cylindrical capacity of $U$ is
$$ {c}^Z(U) := \inf \left \{\pi r^2 \, | \,
\ There \ exists \  \psi \in {\rm Symp}({\mathbb R}^{2n}) \ with \
\psi (U) \subset Z^{2n}(r)  \right  \}.$$
\end{definition}

Note that both the symplectic radius and the cylindrical capacity
satisfy the axioms of Definition~\ref{Def-sym-cap} by the
non-squeezing theorem. Moreover, it follows from
Definition~\ref{Def-sym-cap} that for every symplectic capacity
$c$ and every open set $U \subset {\mathbb R}^{2n}$ we have
$c_B(U) \le c(U) \le c^Z(U)$.

The above axiomatic definition of symplectic capacities is
originally due to Ekeland and Hofer~\cite{EH}. Nowadays, a variety
of symplectic capacities are known to exist. 
For several of the detailed discussions on symplectic capacities we
refer the reader
to~\cite{CHLS},~\cite{Ho},~\cite{HZ},~\cite{L},~\cite{Mc}
and~\cite{V1}.

In this note we mainly concentrate on two important examples of
symplectic capacities which arose from the study of periodic
solutions of Hamiltonian systems. These are the Ekeland-Hofer
capacity $c_{EH}$ introduced in~\cite{EH},~\cite{EH1} and the
Hofer-Zehnder capacity $c_{HZ}$ introduced in~\cite{HZ1}. These
invariants have several applications, among them a new proof of
Gromov's non-squeezing theorem, establishing the existence of closed
characteristics on or near an energy surface, and studying the Hofer
geometry on the group of Hamiltonian diffeomorphisms (see
e.g~\cite{HZ}). Moreover, it is known that on the class of convex
bodies in ${\mathbb R}^{2n}$, these two capacities coincide, and 
can be represented by the minimal symplectic area of a closed
characteristic on the boundary of the convex domain. 
Since in this note we are concerned only with convex sets, we omit
the general definitions of these two capacities, and give a
definition which coincides with the standard ones on the class of
convex domains. This is done in Theorem~\ref{Cap_on_covex_sets}
below. Next we turn to some background on Hamiltonian dynamics.

\subsection{Hamiltonian Dynamics on Convex Domains} \label{sec-Ham-dynm}

Let $U$ be a bounded, connected, open set in ${\mathbb R}^{2n}$ with
smooth boundary containing the origin. A nonnegative function $F :
{\mathbb R}^{2n} \rightarrow {\mathbb R}$ is said to be a {\it
defining function} for $U$ if it satisfies that $\partial U =
F^{-1}(1)$, that $U = F^{-1}([0,1])$, and that $1$  is a regular
value of $F$. Next, let $F$ be a defining function for $U$, and
denote by $X_F = J \nabla F$ the corresponding Hamiltonian vector
field  defined by $i_{X_F} \omega = -dF$. Note that $X_F$ is always
tangent to $\partial U$ since $dF(x) \cdot X_F(x) = - \omega
(X_F(x),X_F(x)) =0$ for all $x \in \partial U$, and hence it defines
a non-vanishing vector field on $\partial U$. It is well known (see
e.g~\cite{HZ}) that the orbits of this vector
field, that is, 
the solutions of the classical Hamiltonian equation $\dot{x} =
X_F(x)$, do not depend, up to parametrization, on the choice of the
Hamiltonian function $F$ representing $\partial U$. Indeed, if $H$
is another defining function for $\partial U$ i.e.,
$$\partial U  = \{ x \ ; \ H(x) = 1 \} = \{x \ ; \ F(x)=1 \} \
{\rm with} \ dH,dF \neq 0 \ {\rm on} \ \partial U,$$ and where $1$
is a regular value of both $F$ and $H$,  then $dF(x) = \lambda(x)
dH(x)$ at every point $x \in
\partial U$, with $\lambda(x) \neq 0$, and therefore
 $X_F = \lambda X_H$  on $\partial U$ where $\lambda \neq 0$.
Thus, the two vector fields have, up to reparametrization, the same
solutions on ${\partial U}$.

The images of the periodic  solutions of the above mentioned
Hamiltonian equation are called the ``closed characteristics'' of
$\partial U$ (where periodic means $T$-periodic for some positive
$T$).
The breakthrough in the global existence of closed characteristics
was achieved by Weinstein~\cite{W} and Rabinowitz~\cite{R} who
established the existence of a closed characteristic on every
convex (and in fact also on every star-shaped) hypersurface in
${\mathbb R}^{2n}$.

We recall the following definition. The action of a $T$-periodic
solution $l(t)$ 
is defined by (see e.g.~\cite{HZ} Page 7):
\begin{equation}\label{action!}
{\cal A}(l) = \int_l \lambda = {\frac 1 2} \int_0^T \langle -J
\dot{l}(t) , l(t) \rangle dt,
\end{equation}
where $\lambda = \sum_1^n x_idy_i$ is the Liouville 1-form whose
differential is $d\lambda = \omega$. This action of a periodic
orbit $l(t)$ is the symplectic area of a disc spanned by the loop
$l(t)$.

In particular, it is a symplectic invariant i.e., ${\cal A}(\psi(l)) = {\cal A}(l)$, for any  $\psi \in {\rm
Symp}({\mathbb R}^{2n})$.

We next introduce the Ekeland-Hofer and the Hofer-Zehnder
capacities, denoted by $c_{EH}$ and $c_{HZ}$ respectively. As stated
above, instead of presenting the general definitions of these two
capacities we present a definition sufficient for our purpose which
coincides with the standard ones on the class of convex domains.
This definition follows from the theorem below, which is a
combination of results from~\cite{EH} and~\cite{HZ}.

\begin{theorem} \label{Cap_on_covex_sets} Let $K \subset {\mathbb
R}^{2n}$ be a convex bounded domain with smooth boundary $\partial
K$. Then there exists at least one closed characteristic $\gamma^*
\subset \partial K$ satisfying
$$ c_{EH}(K) = c_{HZ}(K) = {\cal A}(\gamma^*) =  \min \{ | {\cal A}(l) | \ : \ l \
{\rm is \ a \ closed \ characteristic \ on  \ } \partial K \}$$
\end{theorem}

Such a closed characteristic, which minimizes the action (note that
there might be more than one), is called throughout this text a
``capacity carrier" for $K$. In addition, we refer to the coinciding
Ekeland-Hofer and Hofer-Zehnder capacities on the class of convex
domains as the Ekeland-Hofer-Zehnder capacity, and denote it from
here onwards, when there is no possibility for confusion with a
general capacity, by $c$.

\subsection{Main Results}

A natural question following from the discussion above is whether a
Brunn-Minkowski type inequality holds for the symplectic-size of
sets, which is given by their symplectic capacities. In this paper
we restrict ourselves to the class of convex domains and to the
Ekeland-Hofer-Zehnder capacity. However, we do not exclude the
possibility that the Brunn-Minkowski inequality holds for other
symplectic capacities or other, more general classes of
bodies in ${\mathbb R}^{2n}$. For example, in dimension 2, 
any symplectic capacity agrees with the volume for a large class of
sets in ${\mathbb R}^2$ (see~\cite{Sib}), and hence the
Brunn-Minkowski inequality holds for this class. Also, it is not
difficult to verify that for the linearized ball capacity (for a
definition see~\cite{AMO},~\cite{AO}) the Brunn-Minkowski inequality
holds.

The main result in this paper is the following: Denote by ${\cal
K}^{2n}$ the class of compact convex bodies in ${\mathbb R}^{2n}$
which has non-empty interior.

\begin{theorem}\label{MT}
Let $c$ be the Ekeland-Hofer-Zehnder capacity.
 Then for any $n$, and any
 $K,T \in {\cal K}^{2n}$, one has
\begin{equation} \label{BM-ineq-for-p=1}
c(K+ T)^{\frac 1 2} \geq  c(K)^{\frac 1 2}+ c(T)^{\frac 1 2}.
\end{equation}
Moreover, equality holds if and only if $K$ and $T$ have a pair of homothetic capacity carriers.
\end{theorem}
In fact, Theorem~\ref{MT} is a special case of a slightly more general result
which we now  describe.
For a convex body $K$, denote by $\|x\|_K := \inf\{r: x/r \in K\}$ the
corresponding gauge function. Moreover, we uniquely associate with $K$ its support function
$h_K$ given by:
$$ h_K(u) = \sup \{ \langle x, u \rangle \ : \  x \in K \}, \ \ {\rm
for \ all} \ u \in {\mathbb R}^{2n}.$$ Note that this is no other
than the gauge function of the polar body
\[K^{\circ} = \{ x \in {\mathbb R}^{2n} \ : \ \langle x , y \rangle
\leq 1, \ {\rm for} \ {\rm every} \ y \in K \},\]
or, in the
symmetric case, simply the dual norm $h_K(u) = \| u\|_K^* = \| u
\|_{K^{\circ}}$.

In~\cite{F}, Firey introduced a new operation for convex bodies,
called ``$p$-sum", which depends on a parameter $p \geq 1$ and
extends the classical Minkowski sum. For two convex bodies $K,T
\in {\mathbb R}^{2n}$, both containing the origin, the $p$-sum of
$K$ and $T$, denoted $K+_pT$, is defined via its support function
in the following way:
\begin{equation} \label{p-sum} h_{K+_pT}(u)
= \bigl ( h^p_K(u) + h^p_T(u) \bigr )^{\frac  1 p}, \ \  u \in
{\mathbb R}^{2n}. \end{equation} The convexity of $h_{K+_p T}$
follows easily from Minkowski's inequality. The  case $p = 1$
corresponds to the classical Minkowski sum. Thus, Theorem~\ref{MT}
is a special case of the following:

\begin{theorem}\label{MT_for_p}
Let $c$ be the Ekeland-Hofer-Zehnder capacity.
 Then for any $n$, any $p \geq 1$, and any
 $K,T \in {\cal K}^{2n}$ one has
\begin{equation} \label{BM-ineq-gen}
c(K+_p T)^{\frac p 2} \geq  c(K)^{\frac p 2}+ c(T)^{\frac p 2}.
\end{equation}
Moreover, equality holds if and only if $K$ and $T$ have a pair of homothetic capacity carriers.
\end{theorem}
An interesting corollary of Theorem~\ref{MT} is a symplectic
analogue of the classical isoperimetric inequality comparing volume
and surface area which we now present. For a curve $ \gamma :[0,T]
\rightarrow {\mathbb R}^{2n}$ and a convex body $K$ including $0$ in
its interior we denote by ${\rm length}_{K} (\gamma) = \int_0^T \|
\dot \gamma(t) \|_K dt$ the length of $\gamma$ with respect to the
body $K$. The following corollary is proven in
Section~\ref{Section-Corr&Disscussion}.

\begin{corollary} \label{Iso_ineq}
For any $K,T \in {\cal K}^{2n}$, and $c$ as above,
 \begin{equation} \label{Corollary-1} 4 c(K) c(T) \le ({\rm
length}_{JT^{\circ}} (\gamma))^2, \end{equation} where $\gamma$ is
any capacity carrier of $K$.
\end{corollary}
In Section~\ref{Section-Corr&Disscussion} we explain why Corollary~\ref{Iso_ineq} can be
thought of as a consequence of an isoperimetric-type inequality for
capacities. Note that in the special case where $T$ is the Euclidean unit
ball, Equation~$(\ref{Corollary-1})$ becomes
$$ 4 \pi c(K) \leq  ({\rm length(\gamma)})^2, $$
where $\gamma$ is any capacity carrier for $K$ and where ${\it
length}$ stands for the standard Euclidean length. This last
consequence is known, and can be deduced from the standard
isoperimetric inequality in ${\mathbb R}^{2n}$ combined with the
well known fact that the symplectic area is always less than or
equal to the Euclidean area.

Another special case of Corollary~\ref{Iso_ineq} which can be useful
is the following: let $K$ be a symplectic ellipsoid ${ E}=
\sum_{i=1}^n {\frac {x_i^2 + y_i^2} {r_i^2}},$ where $1 = r_1 \leq
r_2, \ldots \leq r_n$. Equation~$(\ref{Corollary-1})$ implies that
for any $T \in {\cal K}^{2n}$
$$ 4 \pi c(T) \leq ({\rm length}_{JT^{\circ}} (S^1)^2,$$ where
$S^1$ is the  capacity carrier of ${E}$ given by $ x_1^2 + y_1^2 =
1$. Moreover, since the same is true for any symplectic image of
${E}$, we get that
\begin{equation}
 4 \pi c(T) \leq \inf_{\varphi \in {\rm Symp}({\mathbb R}^{2n})}
 ({\rm length}_{JT^{\circ}} (\varphi (S^1)))^2.
\end{equation}
This estimate is sometimes strictly better than other available
estimates, such as volume radius (see~\cite{V},~\cite{AO},~\cite{AMO}).

Next we state another corollary of Theorem \ref{MT}, which improves
 a result previously proved in~\cite{AMO} by other methods. Define
the ``mean-width" of a centrally symmetric convex body $K$ to be
 $$M^*(K) :=  \int_{S^{2n-1}} \max_{y \in K} \langle x
, y \rangle \sigma(dx),$$ where $\sigma$ is the rotationally
invariant probability measure on the unit sphere $S^{2n-1}$. In
other words, we integrate over all unit directions $x$ half of the
distance between two parallel hyperplanes touching $K$ and
perpendicular to the vector $x$. Mean-width is an important
parameter in Asymptotic Geometric Analysis, and is the geometric
version of a central probabilistic parameter, see e.g.
\cite{Pisbook}. We show in Section~\ref{Section-Corr&Disscussion}
below that:

\begin{corollary} \label{Corr-Mean-Width}
For every centrally symmetric convex body $K \subset {\cal K}^{2n}$,
one has
$$ c(K) \leq \pi (M^*(K))
$$ Moreover, equality holds if and only if $K$ is a Euclidean ball.
\end{corollary}

In fact, as the proof will demonstrate, this corollary follows from
standard arguments once we have a Brunn-Minkowski-type inequality.
The same is true for the following result.

\begin{corollary} \label{Corr-intersections}
For any two symmetric convex bodies $K,T \subset {\cal K}^{2n}$, one
has for every $x\in \R^{2n}$ that
\[
c(K\cap (x+T)) \leq c(K \cap T).\] More generally, for any $K, T\in
{cal K}^{2n}$ any $x,y \in \R^{2n}$ and any $0\le \lambda\le 1$, we
have that
\[ \lambda c^{1/2}(K \cap (x+T)) + (1-\lambda)  c^{1/2}(K \cap
(y+T)) \le c^{1/2}(K \cap (\lambda x+(1-\lambda) y+ T)).\]
\end{corollary}

We wish to remark that this note can be considered as a continuation
of the line of work which was presented in \cite{AO} and \cite{AMO},
in which we used methods and intuition coming from the field of
asymptotic geometric analysis to obtain results in symplectic
geometry.

\noindent {\bf Structure of the paper:} The paper is organized as
follows. In Section~\ref{sec_main_ing} we introduce the main
ingredient in the proof of our main theorem. In
Section~\ref{sec-proof-MR} we prove the Brunn-Minkowski inequality
for the Ekeland-Hofer-Zehnder capacity $c$ and characterize the
equality case. In Section~\ref{Section-Corr&Disscussion} we prove
the above mentioned applications of the inequality, and in the last
section we prove a technical lemma.

\noindent {\bf Acknowledgment:} We cordially thank Leonid
Polterovich for very helpful remarks.

\section{The Main Ingredient} \label{sec_main_ing}

In this section we introduce the main ingredient in the proof of
Theorem~\ref{MT_for_p}.

We note that there is no loss of generality in assuming, from here
onwards, that in addition to being compact and with non-empty
interior, all  convex bodies considered also have a smooth boundary
and contain the origin in the interior. Indeed, affine translations
in ${\mathbb R}^{2n}$ are symplectomorphisms, which accounts for the
assumption that the origin is in the interior. Secondly, once we
know the Brunn-Minkowski inequality for smooth convex domains, the
general case follows by standard approximation, since symplectic
capacities are continuous on the class of convex bodies with respect
to the Hausdorff distance (see e.g.~\cite{McS}, Page 376).

The main ingredient in the proof of Theorem~\ref{MT_for_p}, is the
following proposition which is another characterization of the
Ekeland-Hofer-Zehnder capacity, valid for smooth convex sets. Let
$W^{1,p}(S^1,{\mathbb R}^{2n})$ be the Banach space of absolutely
continuous $2 \pi$-periodic functions whose derivatives belong to
$L_p(S^1,{\mathbb R}^{2n})$.

\begin{proposition} \label{main_prop}
Let $c$ be the Ekeland-Hofer-Zehnder capacity. For any convex body $K \subset
{\mathbb R}^{2n}$ with smooth boundary, and any two parameters $p_1
> 1$, $p_2 \geq 1$
\begin{equation}\label{ep}
 c(K)^{\frac {p_2} 2 }= \pi^{p_2} \, \min_{z \in
{\cal E}_{p_1}}
  {\frac 1 {2 \pi}} \int_0^{2 \pi} h_K^{p_2}(\dot{z}(t))dt ,\end{equation}
where $$ {\cal E}_{p_1} = \left \{ z \in W^{1,p_1}(S^1,{\mathbb
R}^{2n}) \ : \  \int_0^{2 \pi} z(t) dt = 0, \  {\frac 1 2} \int_0^{2
\pi} \langle Jz(t), \dot{z}(t) \rangle dt = 1 \right \}.$$
\end{proposition}

In the case where $p_1=p_2=2$, a proof of the above proposition can
be found in~\cite{HZ} and~\cite{MZ}. There the authors use the idea
of dual action principle by Clarke~\cite{C} in order to prove the
existence of a closed characteristic for convex surfaces, a claim
originally due to Rabinowitz~\cite{R} and Weinstein~\cite{W}. For
further discussions on Clarke's dual action principle, and in
particular its use for the proof of existence of closed
characteristics, see e.g.~\cite{Ek} and the references within.

It turns out that the special case $p_1 = p_2 >1$ implies the more
general case of possibly different $p_1>1$, $p_2\ge 1$. That is, we
claim that the following Proposition formally implies Proposition
\ref{main_prop}:

\begin{proposition} \label{special_case_of_main_prop}
For any convex body $K \subset {\mathbb R}^{2n}$ with smooth
boundary, and $p > 1$
\begin{equation}
 c(K)^{\frac {p} 2 }= \pi^{p} \, \min_{z \in
{\cal E}_{p}}
  {\frac 1 {2 \pi}} \int_0^{2 \pi} h_K^{p}(\dot{z}(t))dt.\end{equation}
 \end{proposition}

\noindent{\bf Proof of the implication Proposition
\ref{special_case_of_main_prop} $\Rightarrow$ Proposition
\ref{main_prop}.} Note that for $1<p_2 \le p_1$, one has ${\cal
E}_{p_1} \subset {\cal E}_{p_2}$. Moreover, from H\"{o}lder's
inequality it follows that
\[ \Big(
\frac{1}{2\pi} \int_0^{2\pi} h_K^{p_2}(\dot{z}(t))dt \Big )^{\frac 1 {p_2}}
\le \Big ( \frac{1}{2\pi} \int_0^{2\pi}
h_K^{p_1}(\dot{z}(t))dt \Big )^{\frac 1 {p_1}}.
\]
Therefore, from Proposition~\ref{special_case_of_main_prop} it follows that for $1 < p_2\le p_1$
\begin{eqnarray*}
c(K)^{\frac 1 2} &= &\pi \min_{z\in {\cal E}_{p_1}} \Big(
\frac{1}{2\pi}
\int_0^{2\pi} h_K^{p_1}(\dot{z}(t))dt\Big)^{\frac 1 {p_1} }\\
& \ge & \pi \min_{z\in {\cal E}_{p_1}} \Big( \frac{1}{2\pi}
\int_0^{2\pi} h_K^{p_2}(\dot{z}(t))dt\Big)^{\frac 1 {p_2}}\\
& \ge & \pi \min_{z\in {\cal E}_{p_2}} \Big( \frac{1}{2\pi}
\int_0^{2\pi} h_K^{p_2}(\dot{z}(t))dt\Big)^{\frac 1 {p_2}} =
c(K)^{\frac 1 2}.
\end{eqnarray*}
In particular, we have equality throughout. Similarly,
\begin{eqnarray*}
c(K)^{\frac 1 2} &= &\pi \min_{z\in {\cal E}_{p_1}} \Big(
\frac{1}{2\pi}
\int_0^{2\pi} h_K^{p_1}(\dot{z}(t))dt\Big)^{\frac 1 {p_1} }\\
& \ge & \pi \min_{z\in {\cal E}_{p_2}} \Big( \frac{1}{2\pi}
\int_0^{2\pi} h_K^{p_1}(\dot{z}(t))dt\Big)^{\frac 1 {p_1}}\\
& \ge & \pi \min_{z\in {\cal E}_{p_2}} \Big( \frac{1}{2\pi}
\int_0^{2\pi} h_K^{p_2}(\dot{z}(t))dt\Big)^{\frac 1 {p_2}} =
c(K)^{\frac 1 2}.
\end{eqnarray*}
Thus, we conclude that for any $1<p_1,p_2$
\begin{equation} \label{eq_p_1_p_2} c(K)^{\frac 1 2} = \pi \min_{z\in {\cal E}_{p_1}} \left( \frac{1}{2\pi}
\int_0^{2\pi} h_K^{p_2}(\dot{z}(t))dt\right)^{\frac 1 {p_2}}.
\end{equation}

To complete the proof we need only to explain the case of $p_2=1$.
On the one hand, H\"older's inequality implies that
\[c(K)^{\frac 1 2} = \pi \min_{z\in {\cal E}_{p_1}}
\Big( \frac{1}{2\pi} \int_0^{2\pi}
h_K^{p_1}(\dot{z}(t))dt\Big)^{\frac 1 {p_1} }
 \ge  \pi \min_{z\in {\cal E}_{p_1}}\frac{1}{2\pi}
\int_0^{2\pi} h_K(\dot{z}(t))dt.\]

On the other hand, using that $\lim (\min) \le \min (\lim)$, we can
let $1< p_2 \to 1$ in equation~$(\ref{eq_p_1_p_2})$. By Lebesgue's
dominated convergence theorem we can also insert the limit into the
integral and get that
\[c(K)^{\frac 1 2} \le \pi \min_{z\in {\cal
E}_{p_1}}\frac{1}{2\pi} \int_0^{2\pi} h_K(\dot{z}(t))dt,\]
which completes the proof of Proposition~\ref{main_prop} in the case
of $p_2 = 1$. $\hfill \square$

Before turning to the proof of
Proposition~\ref{special_case_of_main_prop}, which will be our main
objective throughout the rest of this section, let us point out an
important consequence of the above argument which will be helpful
for us later (especially in the proof of the equality case for the
Brunn-Minkowski inequality).

Fix $p_1>1$ and let $\tilde z$ be any path in ${\cal E}_{p_1}$ for
which the minimum is attained in equation~$(\ref{eq_p_1_p_2})$.
Letting $1 \leq p_2 <p_1$ we get that
\begin{eqnarray*}
c(K)^{\frac 1 2} &= &\pi
\Big(\frac{1}{2\pi} \int_0^{2\pi} h_K^{p_1}(\dot{\tilde{z}}(t))dt\Big)^{\frac 1 {p_1} }\\
& \ge & \pi \Big( \frac{1}{2\pi}
\int_0^{2\pi} h_K^{p_2}(\dot{\tilde{z}}(t))dt\Big)^{\frac 1 {p_2}}\\
& \ge & \pi \min_{z\in {\cal E}_{p_1}} \Big( \frac{1}{2\pi}
\int_0^{2\pi} h_K^{p_2}(\dot{z}(t))dt\Big)^{\frac 1 {p_2}} =
c(K)^{\frac 1 2}.
\end{eqnarray*}
In particular, there is equality in the first inequality so that the
$L_{p_1}$ and $L_{p_2}$ norms of the function
$h_K(\dot{\tilde{z}}(t))$ coincide. This clearly implies that this
function is constant in $t$. Another fact which easily follows from
the line of inequalities above is that the minimum is attained on
the {\em same} paths $z$ for all $p_2\ge 1$ (in particular, on a
function which belongs to $\bigcap_{p>1}{\cal E}_p$).

Thus, we have shown that Proposition~\ref{special_case_of_main_prop}
implies the following

\begin{corollary}\label{gulu} {\it
Fix $p_1>1$ and $p_2\ge 1$. Any path $\tilde{z}$ which minimizes
$\int_0^{2 \pi} h_K^{p_2}(\dot{z}(t))dt$ over ${\cal E}_{p_1}$
satisfies that the function $h_K(\dot{z}(t))$ is the constant
function $c(K)/\pi$, and in particular all the $L_p$ norms of the
function $h_K(\dot{z}(t))$ coincide.}
\end{corollary}
After the proof of Proposition~\ref{special_case_of_main_prop}, we
will give a geometrical explanation to this fact, see
Remark~\ref{schnei} below.

We now turn to the proof of Proposition
\ref{special_case_of_main_prop}. We follow closely the arguments,
valid for $p = 2$, in~\cite{HZ} and~\cite{MZ}. Fix $p > 1$ and
consider the functional
\[I_p(z) =
\int_0^{2 \pi} h_K^p(\dot{z}(t))dt\]
defined on the space ${\cal E}_p$, which was defined in the
statement of Proposition \ref{main_prop}.
A key ingredient in the proof is Lemma~\ref{tech_lemma} below, which
we will prove in Section~\ref{prooflemma}, and which gives a
one-to-one correspondence between the so called ``critical points''
of the functional $I_p$ and closed characteristics on $\partial K$.
Before stating the lemma, we must define what we mean by a critical
point of $I_p$, since ${\cal E}_p$ is not closed under all
perturbations.

\begin{definition} \label{def-crit-pt}
An element $z \in {\cal E}_p$ is called a critical point of $I_p$ if
the following holds: For every $\xi \in W^{1,p}(S^1, \R^{2n})$
satisfying $\int_0^{2 \pi} \xi(t)dt = 0$ and $\int_0^{2 \pi} \langle
\xi(t), J \dot{z}(t)\rangle dt = 0$, one has
\[ \int_0^{2 \pi} \langle \nabla h_K^p(\dot{z}(t)),
\dot{\xi} (t) \rangle = 0 \]
\end{definition}

To understand why this definition is natural, first notice that
 the above
condition $\int_0^{2 \pi} \langle \xi(t), J \dot{z}(t)  \rangle dt =
0$ implies that $\int \langle \dot{\xi}(t), J {z} (t)  \rangle dt =
0$ by  integration by parts. Next, consider the element $z_{\eps} =
z + \eps \xi$. It belongs to $W^{1,p}(S^1, \R^{2n})$ and satisfies
the normalization condition $\int_0^{2 \pi} z_{\eps}(t)dt = 0$, but
its action is not normalized to be $1$, thus it is not necessarily
in ${\cal E}_p$. However, its action is close to $1$ with difference
being of order $o(\eps)$. Indeed,
\[ |{\cal A}(z_{\eps}) | = {\frac 1 2} \int_0^{2 \pi}
\langle Jz_{\eps}(t), \dot{z_{\eps}}(t)
\rangle dt =  1 + {\frac {\eps^2} {2}} \int_0^{2 \pi} \langle J
\xi(t), \dot \xi (t) \rangle dt \]
Denote by $z_{\eps}'$ the normalized path:
\[ z_{\eps}' = {\frac {z_{\eps}} {1+ {\frac {\eps^2} {2}} \int \langle J\xi(t) ,
\dot{\xi}(t) \rangle dt}}\] Note that now $z_{\eps}' \in {\cal
E}_p$. For a critical point, it is natural to require that the
difference between $I_p(z)$ and $I_p(z_{\eps}')$ will be of order
$o(\eps)$. Taking the first order approximation we have
\[ I_p(z_{\eps}') =
\int_0^{2 \pi} h_K^p ({\dot{z_{\eps}'}})(t))dt  = \int_0^{2 \pi}
h_K^p(\dot{z}(t))dt  + \eps \int_0^{2 \pi} \langle \nabla
h_K^p(\dot{z}(t)), \dot{\xi} (t) \rangle  + o({\eps^2}),\] and for
the second term on the right hand side to disappear we need exactly
the condition in the definition of a critical point above. In
particular, we emphasize that the minimum of $I_p$ over ${\cal E}_p$ is
attained at a critical point according to our definition, a fact
which will be important in the proof. With the definition in hand,
we may formulate the following lemma, which for $p=2$ appears
in~\cite{HZ}, Pages 26-30. For the sake of completion we include its
proof for general $p>1$ in Section~\ref{prooflemma}.

\begin{lemma} \label{tech_lemma}
Let $K \subset {\mathbb R}^{2n}$ be a convex body with smooth boundary, and fix $p>1$. Each critical point of
the functional $I_p(z)$ satisfies the Euler equation
\begin{equation} \label{Euler_eq} \nabla h_K^p ({\dot {z}}) = {\frac
p 2} \,  \lambda J \,  z + \alpha, \qquad {\rm where} ~~~ \lambda =
I_p(z),
\end{equation}
for some fixed vector $\alpha$ (which may be different for different
critical points), and vice versa: each point $z$ satisfying
Equation~$(\ref{Euler_eq})$ is a critical point of $I_p$. Moreover,
the
 functional $I_p(z)$
achieves its minimum i.e., there is $\tilde z \in {\cal E}_p$ such
that \[ I_p(\tilde z) = \inf_{z \in {\cal E}_p} I(z) = \int_0^{2 \pi}
h_K^p( \dot{\tilde z}(t)) dt= \tilde \lambda \neq 0,\] and in
particular, ${\tilde{z}}$ satisfies Equation~$(\ref{Euler_eq})$.
\end{lemma}

\noindent {\bf Proof of Proposition~\ref{special_case_of_main_prop}} The idea is as
follows: we define an invertible mapping ${\cal F}$ between critical
points $z$ of $I_p(z)$ and  closed characteristics $l$ on the boundary
of $K$. Moreover, we will show that the action of $l = {\cal F}(z)$
is a simple monotone increasing function of $I_p(z)$. In particular,
the critical point $z$ for which the minimum of $I_p(z)$ is attained
is mapped to the closed characteristic minimizing the action. Since
the minimal action of a closed characteristic is exactly the
Ekeland-Hofer-Zehnder capacity, the result will follow.

To define the mapping ${\cal F}$, let $z: S^1\to R^{2n}$ be a
critical point of $I_p$.
 In particular from Lemma
\ref{tech_lemma} we have that
 \begin{equation} \label{Euler_eqt}
\nabla h_K^p ({\dot {z}}) = {\frac p 2} \,  \lambda J \, { z} +
\alpha,
\end{equation}
for some vector $\alpha$ and $\lambda = I_p(z)$. We will use the
Legendre transform in order to define an affine linear image of
$z$ which is a closed characteristic on the boundary $\partial K$
of $K$, which we will then define as ${\cal F}(z)$. Recall that
the Legendre transform is defined as follows: For $f : {\mathbb
R}^{n} \rightarrow {\mathbb R}$, one defines
$$ {\cal L} f(y) = \sup_{x \in {\mathbb R}^n} [ \langle y , x
\rangle - f(x) ], \ \forall y \in {\mathbb R}^n.$$
It is not hard to
check that
\[ ({\cal L} (h_K^p))(v) =
{\frac {p^{1-q}} q} h_{K^\circ}^q (v), \] where $p^{-1} + q^{-1} =1
$
and  $K^{\circ}$
is, as before, the polar body of $K$. Note that $h^q_{K^{\circ}}$ is
a defining function of $K$ (that is, $K$ is its 1-level set) which
is homogeneous of degree $q$. After applying the Legendre transform
and using the fact that $v =  \nabla h_K^p(u)$ is inverted
point-wise by $u = \nabla {\cal L} h_K^p(v)$ equation
$(\ref{Euler_eq})$ becomes:
\[ {\dot {z}} =   {\frac  {p^{1-q}}  q} \nabla h_{K^{\circ}}^q \Bigl
({\frac p 2} \,  \lambda J \, \tilde z + \alpha \Bigr)  =  \nabla
h_{K^{\circ}}^q \Bigl ({\frac {q^{ {\frac 1 {1-q}}  }} 2} \, \lambda
J \,  z + {\frac {\alpha q^{ {\frac 1 {1-q}}}} p} \Bigr). \]
 Next, let
\begin{equation} \label{eq_20} l = {\kappa} \Bigl ({\frac
{q^{ {\frac 1 {1-q}}}} 2} \,  \lambda J \, z + {\frac {\alpha q^{
{\frac 1 {1-q}}}} p} \Bigr),
\end{equation}
where $\kappa$ is a positive normalization constant which we will
readily choose. Differentiating~$(\ref{eq_20})$ we see  that $l$
satisfies the following Hamiltonian equation. \begin{equation}
\label{l-Ham-eq} {\dot {l}} = {\frac \kappa 2}\, {{q^{ {\frac 1
{1-q}}}} } \,\lambda \, J \nabla h_{K^{\circ}}^q ( { {l}/ \kappa})
= {\frac {\kappa^{2-q}} 2} \, {{q^{ {\frac 1 {1-q}}}} } \, \lambda
J \nabla h_{K^{\circ}}^q (
 l).\end{equation}
 Note that $l$ is a periodic trajectory of the
Hamiltonian equation corresponding to the Hamiltonian function
$h_{K^{\circ}}^q$. Since we ask $l \in {\partial K}$ we need to
choose $\kappa$ such that $l$ will lie in the energy level
$h_{K^{\circ}}^q=1$. For this purpose, note that since
$h_{K^{\circ}}^q$ is homogeneous of degree $q$ we obtain from
Euler's formula that
\begin{eqnarray*} {\frac 1 {2 \pi}} \int_0^{2\pi} h_{K^{\circ}}^q(
 l(t))dt  & = & {\frac 1 { 2 \pi q} } \int_0^{2\pi} \langle \nabla
h_{K^{\circ}}^q(  l(t)) ,  l(t) \rangle dt = -{\frac {\kappa^{q-2}
q^{\frac 1 {q-1}} } { \pi  \lambda q} } \int_0^{2\pi} \langle J
\dot{  l}(t) , l(t) \rangle dt
\\ & = & {\frac { {{q^{ {\frac 1 {1-q}} }}} \kappa^{q}  } { 4 \pi q  \lambda}}
\int_0^{2\pi} \langle   \lambda {\dot {{z}}}(t) , \lambda J z(t) +
{\frac {2 \alpha}  p} \rangle dt = {\frac { \kappa^{q}
 \lambda q^{\frac q {1-q}}} {2 \pi   }}
\end{eqnarray*}
which is equal to $1$ if we choose $\kappa = ({ {2 \pi } / {
\lambda}})^{\frac 1 q} q^{\frac 1 {q-1}} $. Therefore, for this
value of $\kappa$ we have that
\begin{equation} \label{new_eq_for_l}  l
= \Big( {\frac {2 \pi} \lambda} \Big )^{\frac 1 q} \Big( {\frac
\lambda 2} Jz + {\frac {\alpha} p} \Big), \end{equation} is a closed
trajectory of the Hamiltonian equation  corresponding to the
Hamiltonian $h_{K^{\circ}}^q$ on $\partial K$. This $l$ we denote by
${\cal F}(z)$. (To be completely formal, to agree with the way
closed characteristics were defined, we let ${\cal F}(z)$ be the
image of $l$ in $\R^{2n}$.) Below we will show that this mapping is
invertible, and compute ${\cal F}^{-1}$.

Next we derive the relation between ${\cal A}(l)$ and $\lambda =
I_p(z)$.
 Using Euler's formula again, and the above value of $\kappa$ we
 conclude that
 $$ {\cal A} (l) = {\frac 1 2} \int_0^{2 \pi} \langle -J {\dot l}(t) ,l(t) \rangle dt =
  {\frac {\kappa^2} 8}  q^{\frac 2 {1-q}} {\lambda}^2  \int_0^{2 \pi} \langle
  {\dot z}(t) ,  J  z(t) + {\frac {2 \alpha} {\lambda p}} \rangle dt
 = {4}^{-{\frac {1} p}}
  \left ( {\pi} \right )^{\frac 2 q}  \lambda^{\frac 2 p} $$
 Equivalently,
\begin{equation}\label{lambdaval}
{\cal A}^{\frac p 2} (l) =  1/2 (\pi)^{\frac p q} \lambda.
\end{equation}

In order to show that the map  ${\cal F}$ (we should actually write
${\cal F}_{p}$ as it depends on $p$, but we omit this index so as
not to overload notation) is indeed one-to-one and onto, we now
define ${\cal F}^{-1}$. Starting now with a closed characteristic
$\Gamma$ on $\partial K$, it is not difficult to check that we may
assume using a standard re-parametrization argument that it is the
image of a loop $l$ with $l:[0,2 \pi] \rightarrow {\mathbb R}^{2n}$
and $\dot l = d J \nabla h_{K^{\circ}}^q(l)$, for some constant $d$.
Next, we define
\[{\cal F}^{-1}(l) = J^{-1} \Bigl ((\pi d q)^{-1/2} \Bigl (l - {\frac
1 {2\pi}} \int_0^{2\pi} l(t)dt \Bigr ) \Bigr ).\] We will show that
this map is mapping closed characteristics to critical points. Set
$z={\cal F}^{-1}(l)$. It is easy to check that $\int_0^{2 \pi}
z(t)dt=0$. The fact that $z \in W^{1,p}(S^1, {\mathbb R}^{2n})$
follows from the boundedness of $l$ (as ${\rm Image}(l) \in \partial
K$ is bounded) and the following argument: since $\dot z = {C_1}
\cdot \nabla h_{K^{\circ}}^q(l)$ for some constant $C_1$ and,
$\nabla h_{K^{\circ}}^q$, being homogenous of degree $q-1$, satisfy
$|\nabla h_{K^{\circ}}^q(x)| \le C_2 |x|^{q-1}$ for some constant
$C_2$, we conclude that for some constants $C_3$ and $C_4$,
\[
\int_0^{2\pi} |\dot{z}(t)|^pdt = C_3 \int_0^{2\pi} |\nabla
h_{K^{\circ}}^q(l(t))|^{\frac q {q-1}}dt \le C_4 \int_0^{2\pi}
|l(t)|^qdt < \infty.
\]
Moreover
\begin{eqnarray*} {\cal A}(z) & = & {\frac 1 2} \int_0^{2\pi}
\langle \dot z(t), Jz(t) \rangle dt = {\frac 1 2} (\pi d q)^{-1}
\int_0^{2\pi} \langle J^{-1} \dot l(t), l(t) \rangle dt \\ & = &
{\frac 1 2} (\pi d q)^{-1} \int_0^{2\pi} \langle d \nabla
h_{K^{\circ}}^q(l(t)), l(t) \rangle dt = {\frac 1 {2 \pi}}
\int_0^{2\pi} h_{K^{\circ}}^q(l(t)) =1,
\end{eqnarray*}
where the next to last inequality follows from Euler's formula.
Finally, note that
\[ \dot z = (\pi d q)^{-\frac {1} 2} J^{-1} \dot l = (\pi
d q)^{-\frac {1} 2} d \nabla h_{K^{\circ}}^q(l)= (\pi d q)^{-\frac
{1} 2} d \nabla h_{K^{\circ}}^q \Bigl((\pi d q)^{\frac 1
2}Jz+\int_0^{2\pi} l(t)dt \Bigr)\]
Using the Legendre transform as before we get that
$$ \nabla h_K^p( \dot z) = \alpha Jz + \beta,$$ where
$\alpha$ and $\beta$ are constants (depending on $d$ and $q$).
Moreover,
$$ I_p(z) = \int_0^{2 \pi} h_K^p(\dot{z}(t))dt = {\frac 1 p}
\int_0^{2 \pi} \langle \nabla h_K^p(\dot{z}(t)), \dot z(t) \rangle
dt = {\frac {\alpha} p} \int_0^{2 \pi} \langle  Jz(t), \dot z(t)
\rangle dt ={\frac {2 \alpha} p}
$$
From Lemma~\ref{tech_lemma} it now follows that $z$ is a critical
point of $I_p$ and hence the map ${\cal F}^{-1}$ is well defined.
It is not difficult to show that for every critical point, ${\cal
F}^{-1} {\cal F} (z) = z$ and that for every closed characteristic,
${\cal F} {\cal F}^{-1} (l) = l$, and we omit this computation.

This one-to-one correspondence, and the monotone relation between
${\cal A}({\cal F}(z))$ and $I_p(z)$, implies that the for $\tilde z$,
a critical point for which the minimum of $I(z)$ is attained, its
``dual'' $\tilde l = {\cal F}(\tilde{z})  $ has minimal
 action among all closed characteristics $l$.
This fact together with Theorem~\ref{Cap_on_covex_sets} (which we
consider, for the purpose of this note, as the definition of
Ekeland-Hofer-Zehnder capacity of convex domains) implies that
\begin{equation}\label{givemeaname}
c(K)^{\frac p 2} = {\cal A}^{\frac p 2} (\tilde l) =
 ({1/2})^{p } (2 \pi)^{\frac p q} \lambda = \pi^p {\frac 1 {2 \pi}}
 \int_0^{2 \pi} h_K^p( {\dot {\tilde z}}(t))dt. \end{equation}
The proof of the proposition is now complete. $\hfill \square$

\begin{remark}\label{theFremark} {\rm
In the proof we have shown the following fact, which we will use
later in the note once again: For every $p$, there exists an
invertible mapping ${\cal F} ( = {\cal F}_{p})$, mapping critical
points of $I$ (which by Lemma~\ref{tech_lemma} are exactly the loops
satisfying equation~$(\ref{Euler_eq})$), to closed characteristics
on $\partial K$, and moreover, satisfying
\begin{equation}
{\cal A}^{\frac p 2} ({\cal F}(z)) = ({1/2})^{p } (2 \pi)^{\frac p
q} I_p(z).
\end{equation} }
\end{remark}

\begin{remark}\label{schnei} {\rm
The fact exhibited in Corollary~\ref{gulu}, which might seem
surprising at first, can be geometrically explained (for the case
$p>1$) after the above construction has been made. Indeed, recall
that (fixing some $p>1$) for every path $z \in {\cal E}_p$ which is
a critical point of $I_p$ there corresponds a path $l \in \partial
K$, $l = {\cal F}_p(z)$, which is a linear image of $z$. Moreover,
by the above formulae, putting all the constants together and naming
them $A_1,A_2,A_3$, we have (using that $h_{K^{\circ}}(l) = 1$)
\[\dot z = A_1 J^{-1}\dot{l} = A_2 \nabla h_{K^{\circ}}^q (l) =
A_3  \nabla h_{K^{\circ}} (l).\]
Thus,  we have $h_K(\dot{z}) = A_3h_K( \nabla h_{K^{\circ}} (l))$.
However, a simple fact from Convexity (See \cite{Sch}, Corollary
1.7.3. Page 40) is that for a convex body $T$ and a vector $0\neq
u\in \R^n$, the gradient of the dual norm in a certain direction is
exactly the support of the body in this direction. More formally:
\[ \nabla h_{T}(u)= \{ x \in T : h_T(u) = \langle x,
u\rangle\},\] which in our case, by smoothness, is simply one point
(on the boundary of $T$, of course).
Thus, in particular, $h_K(\nabla h_{K^{\circ}}(u)) =
h_{K^{\circ}}(u)$,
and we see that
\[h_K(\dot{z}) = A_3 h_{K^{\circ}}(l) = A_3,\]
is constant and does not depend on $t$, as claimed in
Corollary~\ref{gulu} and proven there by other means. }
\end{remark}

\section{Proof of the Main result} \label{sec-proof-MR}

In this section we use Proposition \ref{main_prop} to prove our main
theorem.

\begin{proof}[{\bf Proof of Theorem~\ref{MT_for_p}.}]
Fix $p_1 > 1$. It follows from equation~$(\ref{p-sum})$
that for every $z \in {\cal
E}_{p_1}$ and $p \geq 1$ \begin{equation}  \label{eq1-prf-of-MT}
{\frac 1 {2 \pi}}  \int_0^{ 2 \pi}
h^{p}_{K+_{p} T}(\dot{z}(t)) dt = {\frac 1 {2 \pi}}  \int_0^{2 \pi}
\ h^{p}_K(\dot{z}(t) )dt
+ {\frac 1 {2 \pi}}  \int_0^{2 \pi} h^{p}_T(\dot{z}(t))dt \end{equation}
By multiplying
both sides of the above equation by $\pi^{p}$,
taking the minimum over all $z \in {\cal E}_{p_1}$, 
and applying Proposition~\ref{main_prop} above, we conclude that for
every $p \geq 1$
\begin{equation} \label{eq1}  c (K+_{p} T)^{\frac {p} 2}
\geq c(K)^{\frac {p} 2} +  c(T)^{\frac  {p} 2}.
\end{equation}
In particular, for $p=1$ we get
the Brunn-Minkowski inequality for the Ekeland-Hofer-Zehnder capacity.

We turn now to prove the equality case. We start by proving that if
$K$ and $T$ have homothetic capacity carriers, then equality holds
in~$(\ref{eq1})$ for every $p\geq1$.

Let $\Gamma_K \subset
\partial K$ be a capacity carrier for $K$.
As in the proof of Proposition~\ref{special_case_of_main_prop}
above, we can choose a parameterized curve representing $\Gamma_K$
via $\Gamma_K = {\rm Image} \ l_K$, where $ {\dot {l}_K} = d_K J
\nabla h_{K^{\circ}}^2(l_K)$ and $l_K(0) = l_K(2 \pi)$. Moreover, it
follows from the proof of
Proposition~\ref{special_case_of_main_prop} (say in the case $p=2$)
that for every such $l_K$ there is corresponding minimizer $z_k \in
{\cal E}_2$ of the functional $I_2$:
$$z_K
:= {\cal F}^{-1}(l_K) = J^{-1} \bigl ((2 \pi d_K)^{-1/2} \bigl (l_K
- {\frac 1 {2\pi}} \int_0^{2\pi} l_K(t)dt \bigr ) \bigr ),$$ such
that $$ c(K)^{\frac 1 2} = \pi \min_{z \in {\cal E}_2} \Big ( {\frac
1 {2\pi}}\int_0^{2\pi} h_{K}^2(\dot z(t))dt \Big )^{\frac 1 2} = \pi
\Big( {\frac 1 {2\pi}} \int_0^{2 \pi} h_{K}^2(\dot z_K(t))dt \Big
)^{\frac 1 2}$$ Moreover, combining this with
Proposition~\ref{main_prop} and Corollary~\ref{gulu} we conclude
that for every $p\geq1$:
 $$ c(K)^{\frac p 2} = \pi^p \min_{z \in {\cal E}_2}
 {\frac 1 {2\pi}}\int_0^{2\pi} h_{K}^p(\dot z(t))dt  = \pi^p
 {\frac 1 {2\pi}} \int_0^{2 \pi} h_{K}^p(\dot z_K(t))dt.$$
Similarly, let $\Gamma_T$ be a capacity carrier for $T$, set $l_T$
the corresponding parameterized curve which represents $\Gamma_T$,
and let $z_T = {\cal F}^{-1}(l_T) \in {\cal E}_2 $ be the
corresponding critical point of $I_2$ which satisfies
$$ c(T)^{\frac p 2} = \pi^p \min_{z \in {\cal E}_2}
{\frac 1 {2\pi}}\int_0^{2\pi} h_{T}^p(\dot z(t))dt  = \pi^p
 {\frac 1 {2\pi}} \int_0^{2 \pi} h_{T}^p(\dot z_T(t))dt.$$
Note that in order to have equality in~$(\ref{eq1})$, it is enough
to show that $z_K = z_T$. To this end we observe that since
$\Gamma_K$ and $\Gamma_T$ are homothetic, there exist two constants
${\alpha}$ and $\beta$ such that $l_T = \alpha l_K + \beta$. This
implies that $z_T = ({\frac {d_K} {d_T}})^{1/2} \alpha z_K$ and that
${\cal A}(l_T) = \alpha^2 {\cal A}(l_K)$.
 Moreover, since $\Gamma_K$ and $\Gamma_T$ are capacity carriers
of $K$ and $T$ respectively, it follows that ${\cal A}(l_K) = c(K)$
and ${\cal A}(l_T) = c(T)$. Hence, we conclude that $\alpha = (
{\frac {c(T)} {c(K)}})^{1/2}$. On the other hand
\begin{eqnarray*} c(K) & = &{\cal A}(l_K) = {\frac 1 {2\pi}}
\int_0^{2\pi} \langle l_K(t), J l_K(t) \rangle dt = {\frac 1 {2\pi}}
\int_0^{2\pi} \langle l_K(t), d_K \nabla h_{K^{\circ}}^2(l_K(t))
\rangle dt \\ & = & {\frac {2 d_K} {2\pi}} \int_0^{2\pi}
h_{K^{\circ}}^2(l_K(t)) dt = 2 d_K,
\end{eqnarray*} and similarly $c(T) = 2 d_T$. This implies that
$z_K = z_T$ and hence we have an equality in~$(\ref{eq1})$ for every
$p \geq 1$ as required.

Next, we assume that equality holds in~$(\ref{eq1})$ for some $p
\geq 1$ and prove that $K$ and $T$ have homothetic capacity
carriers.

Let $p_1 > 1$ and $p \geq 1$. Note that equality in~$(\ref{eq1})$
implies that
$$ \min_{z \in {\cal E}_{p_1}} \int_0^{2\pi} h_{K+_pT}^p(\dot
z(t))dt = \min_{z \in {\cal E}_{p_1}} \int_0^{2\pi} h_{K}^p(\dot
z(t))dt + \min_{z \in {\cal E}_{p_1}} \int_0^{2\pi} h_{T}^p(\dot
z(t))dt$$ This in turn implies that there exists ${\tilde z} \in
{\cal E}_{p_1} $ such that
$$ \min_{z \in {\cal E}_{p_1}} \int_0^{2\pi} h_{K}^p(\dot z(t))dt = \int_0^{2\pi} h_{K}^p(\dot {\tilde z}(t))dt \ \ {\rm and} \ \
\min_{z \in {\cal E}_{p_1}} \int_0^{2\pi} h_{T}^p(\dot z(t))dt =
\int_0^{2\pi} h_{T}^p(\dot {\tilde z}(t))dt$$
 Combining this fact with Proposition~\ref{main_prop},
 we conclude that
$$
\min_{z \in {\cal E}_{p_1}} \Big ( \int_0^{2\pi} h_{K}^{p_1}(\dot
z(t))dt \Big )^{\frac 1 {p_1}}  = \min_{z \in {\cal E}_{p_1}} \Big (
\int_0^{2\pi} h_{K}^p(\dot z(t))dt \Big )^{\frac 1 p}  = \Big (
\int_0^{2\pi} h_{K}^p(\dot {\tilde z}(t)) \Big )^{\frac 1 p}
 $$
In other words, ${\tilde z} \in {\cal E}_{p_1}$ is a critical point
of the functional $I^K_{p_1}=\int_0^{2\pi} h_{K}^{p_1}(\dot z(t))dt
$
 defined on the space ${\cal E}_{p_1}$, where $p_1 > 1$. It follows from the proof of
 Proposition~\ref{special_case_of_main_prop} together with~$(\ref{new_eq_for_l})$ and~$(\ref{lambdaval})$, that
for such ${\tilde z}$ their is a
 corresponding ${\tilde l}_K$
which satisfies
$$ {\tilde l}_K = \Bigl ({\frac {2 \pi} {\lambda_K}}\Bigr)^{\frac 1 {q_1}} \Bigl ( {\frac \lambda 2} J {\tilde z} + {\frac {\alpha_K} {p_1}} \Bigr) = c(K)^{\frac 1 2}J {\tilde z} + A_K,$$
where $A_K$ is a constant which depends on $K$, $p_1$,
and $q_1^{-1} + p_1^{-1} = 1$. Similarly, since ${\tilde z}$ is also
a critical point of $I^T_{p_1}$, we have that ${\tilde l}_T =
c(T)^{\frac 1 2}J {\tilde z} + A_T$. We conclude that ${\tilde l}_T
= \alpha {\tilde l}_K + \beta$ where $\alpha = c(T)^{\frac 1
2}/c(K)^{\frac 1 2}$. This implies that
 $K$ and $T$ have homothetic capacity carriers and the
 proof of Theorem~\ref{MT_for_p} is now complete.
\end{proof}

\section{Corollaries of the Main Theorem}\label{Section-Corr&Disscussion}

In this section we prove
Corollaries~\ref{Iso_ineq},~\ref{Corr-Mean-Width}
and~\ref{Corr-intersections}. We start with the proof of
Corollary~\ref{Iso_ineq}. As in the case of the classical
isoperimetric inequality, which connects the surface area of a body
and its volume, the Brunn-Minkowski inequality is useful in
obtaining a lower bound for the derivative of the volume-type
function. This follows directly from Theorem~\ref{MT} above and the
following computation: for any convex bodies $K$ and $T$ in
${\mathbb R}^{2n}$ and any $\eps>0$,
\begin{eqnarray} \label{iso-inq-1}
\frac{c (K+\eps T)^{\frac 1 2} - c(K)^{\frac 1 2}}{\eps} &\ge&
\frac{c(K)^{\frac 1 2}+\eps c(T)^{\frac 1 2} - c(K)^{\frac 1
2}}{\eps} = c(T)^{\frac 1 2}.
\end{eqnarray}
The limit on the left hand side as $\eps \to 0^+$ can be thought of
as a ``directional derivative'' of $c^{1/2}$ in the ``direction''
$T$. Note that this argument holds for any symplectic capacity for
which one is able to show that the Brunn-Minkowski inequality holds.
However, to get a meaningful result, one must find a geometric
interpretation for the so-called derivative which one arrives at.
To be more precise, let us define, for a convex body $T$, the
functional $d_T(K)$ by
\[ d_T(K) = \lim_{\varepsilon \rightarrow
0^+} \frac{c (K+\eps T) - c(K)}{\eps}.\]
Inequality~$(\ref{iso-inq-1})$ implies the following easy corollary:
\begin{corollary}\label{dee} For every convex body $K \subset {\mathbb R}^{2n}$, one has
\[ d_T(K) = 2c(K)^{\frac 1 2} {\frac d {d \varepsilon}}
c(K+\varepsilon T)^{\frac 1 2}|_{\eps = 0^+} \geq 2 c(K)^{\frac 1
2}c(T)^{\frac 1 2}.\]
\end{corollary}
The only part which requires justification is the existence of the
limit. Let us show that $c^{1/2}(K+\eps T)$ has derivative at $\eps
= 0^+$ (which is, since $c(K)\neq 0$, the same as showing that
$c(K+\eps T)$ has a derivative): Let $s<t$, note that
\[ \frac{c(K+ tT)^{\frac 1 2} - {c(K)^{\frac 1 2}}}{t} \le
\frac{c(K+ sT)^{\frac 1 2} - {c(K)^{\frac 1 2}}}{s},
\]
is equivalent to
 \[ (s/t)(c(K+ tT)^{\frac 1 2}) + (1 - s/t)c
(K)^{\frac 1 2} \le  c(K+ sT)^{\frac 1 2}, \]
which follows from the Brunn-Minkowski inequality. Hence, the
expression in the limit is a decreasing function of $\eps>0$, and
converges to its supremum as $\eps \to 0^+$, provided it is bounded.
Showing that it is bounded is simple, since $T\subset R K$ for some
$R>0$ (which can be huge, and may depend on the dimension) and thus
\[ \frac{c(K+ tT)^{\frac 1 2} - {c(K)^{\frac 1 2}}}{t}
\le  \frac{c (K+ tRK)^{\frac 1 2} - {c (K)^{\frac 1 2}}}{t}  = R.\] This
completes the proof of Corollary \ref{dee}. $\hfill \square$

One way to use Corollary~\ref{dee} is to find a geometric
interpretation of the derivative of the capacity, $d$. Roughly
speaking, if $c$ is a symplectic ``volume'', $d$ should be a kind of
symplectic ``surface-area''. Instead, what we do below is to bound
$d$ from above by an expression with a clear geometric meaning:
minimal length of loops in a certain norm (as in the statement of
Corollary~\ref{Iso_ineq}), and then Corollary~\ref{dee} gives a
lower bound, in terms of capacity, of this expression.

We fix $\eps > 0$, $p_1 > 1$ and $p_2=1$, and denote by $\tilde{z}
\in {\cal E}_{p_1}$ any path on which the minimum in
Equation~$(\ref{ep})$ is attained for $c(K)^{\frac 1 2}$. We
compute:
\begin{eqnarray*}
c(K+ \eps T)^{\frac 1 2} &=&   \pi \, \min_{z \in {\cal E}_{p_1}}
{\frac 1 {2\pi}} \int_0^{2 \pi} h_{K+\eps T}(\dot{z}(t))dt \\
&=& \pi \, \min_{z \in {\cal E}_{p_1}} {\frac 1 {2\pi}} \int_0^{2
\pi} h_K(\dot{z}(t))dt + \eps h_T(\dot{z}(t))dt  \\
&\le & \pi \, {\frac 1 {2\pi}} \int_0^{2
\pi} h_K(\dot{\tilde{z}}(t))dt + \eps h_T(\dot{\tilde{z}}(t))dt  \\
&= & c(K)^{\frac 1 2} + {\frac \eps {2}} \int_0^{2 \pi}
h_T(\dot{\tilde{z}}(t))dt.
\end{eqnarray*}
Rearranging (for fixed $\eps>0$ and $p>1$), and applying
equation~$(\ref{BM-ineq-for-p=1})$, one gets that for any such $\tilde{z}$
\[ c(T)^{\frac 1 2} \leq  {\frac{c(K+\eps T)^{\frac 1 2}-c(K)^{\frac 1
2}}{\eps}} \le {\frac 1 2} \int_0^{2 \pi}
h_T(\dot{\tilde{z}}(t))dt.\]
Note that the middle expression is a decreasing function of $\eps$,
which as $\eps\to \infty$ converges to the left hand side, and as
$\eps \to 0^+$ converges to $d_T(K)/2\sqrt{c(T)}$.

Next, note that ${\tilde z}$ is a critical point of the functional
$I_{p_1}^K(z) = \int_0^{2\pi} h_K^{p_1}({\dot z}(t))dt$ defined on
the space ${\cal E}_{p_1}$ as well (see Corollary \ref{gulu} above
and the reasoning before it). Hence, we can use the transformation
${\cal F}$ defined in the proof of
Proposition~\ref{special_case_of_main_prop}, to map ${\tilde z}$ to
the corresponding
 capacity carrier $\tilde l$ of $K$. Moreover, from
 equalities~$(\ref{new_eq_for_l})$
 and~$(\ref{lambdaval})$ it follows that
$$ {\dot {\tilde l}} = c(K)^{\frac 1 2} J {\dot {\tilde z}}.$$

Thus, the above inequality takes the form
 \begin{eqnarray*}  c(T)^{\frac 1 2} \leq
\frac{c(K+ \eps T)^{\frac 1 2} -  c(K)^{\frac 1 2}}{\eps} \le {\frac
{1} {2}} c(K)^{-\frac {1} 2}
\int_0^{2 \pi} h_T(J^{-1}\dot{\tilde{l}}(t))dt.
\end{eqnarray*}

Since there is a one-to-one correspondence between critical points
${\tilde z}$ of the functional $I_{p_1}^K$ and closed
characteristics $\tilde l$ on $\partial K$,
 we may in fact write the above inequality as
 \begin{eqnarray} \label{eq_blablabla} c(T)^{\frac 1 2} \leq
\frac{c(K+ \eps T)^{\frac 1 2} -  c(K)^{\frac 1 2}}{\eps} \le {\frac
{1} {2}} c(K)^{-\frac {1} 2} \inf_{\tilde{l}}
\int_0^{2 \pi} h_T(J^{-1} \dot{\tilde{l}}(t))dt, 
\end{eqnarray}
where the infimum runs over all the loops $l$ which are images under
${\cal F}$ of $\tilde{z}$ minimizing equation~$(\ref{ep})$ for
$c(K)^{\frac 1 2}$ i.e., all the capacity carriers of $K$.

Thus, we arrive at
\[ 4 c(K) c(T) \le ({\rm length}_{JT^{\circ}} (l))^2,\]
for any capacity carrier $l$ on $\partial K$, proving
Corollary~\ref{Iso_ineq}. We may also take the limit in~$(\ref{eq_blablabla})$ as $\eps \to 0^+$
to see that

\begin{corollary}
For any $n$, any $K,T \in {\cal K}^{2n}$, and any capacity carrier
on $\partial K$, we have that
\begin{eqnarray*}
d_T(K) \le
 \inf_{\tilde{l}}  \int_0^{2 \pi} h_{T}(J^{-1}
\dot{\tilde{l}}(t))dt =  {\rm length}_{JT^{\circ}} (l).
\end{eqnarray*}
\end{corollary}

Next we turn to the relation between the capacity and the Mean-Width
of a body.

\noindent {\bf Proof of Corollary~\ref{Corr-Mean-Width}.} We denote
by $U(n)$ the group of unitary transformations in ${\mathbb C}^n
\simeq {\mathbb R}^{2n}$. Note that $c(UK) = c(K)$ for any unitary
operator $U \in U(n)$. The Brunn-Minkowski inequality thus implies
that for $U_1, U_2 \in U(n)$
\[ c(K) \leq c \left (\frac{U_1K + U_2K}{2} \right ) ,\]
and by induction also
\[  c(K) \leq c \left ( \frac{1}{N} \sum_{i=1}^N U_iK \right ) .\]
Further, this is true also if we integrate (with respect to
Minkowski addition) along the unitary group with respect to the
uniform Haar measure $d \mu$ on this group:
\[ c(K) \leq c \left ( \int UK d \mu (U) \right ).\]
However, it is not hard to see that the integral on the left hand
side is simply a Euclidean ball of some radius, since it is
invariant under rotations $U \in U(n)$, and further, it is easy to
determine its radius since $M^*(K) = M^*(UK)$, for $U \in U(n)$, and
$M^*$ is an additive function with respect to Minkowski addition
(for more details see~\cite{Pisbook}). Thus we have $\int UK d \mu
(U) = M^*(K) B_2^{2n}$, where $B_2^{2n}$ is the Euclidean unit ball,
and the inequality above translates to
\begin{equation}  \label{eq_m_star} c(K) \leq c( M^*(K)B_2^{2n})={\pi}(M^*(K))^2. \end{equation}

Next, we turn to the characterization of the equality case. Let $L$
denote the family of all capacity carriers on $\partial K$. Since we
assume equality between the left and right hand side, we get
equality throughout the following
\[ c(K) \le \frac{1}{{N^2}}c(U_1K + \cdots + U_NK) \le
c(M^*B_2^{2n}),\] for any $N$ and $U_1, \ldots, U_N \in U(n)$.

Applying the same argument as in the proof of Brunn-Minkowski for
two summands (Theorem~\ref{MT}), this time to $N$ summands $K_1,
\ldots K_N$ for arbitrary $N$, we get that the equality conditions
becomes the following: there exists $N$ homothetic capacity carriers
$l_i \subset \partial K_i$. Since in this case $K_i = U_iK$, we know
that capacity carriers on $K_i$ are images by $U_i$ of capacity
carriers on $K$. Moreover, since $K$ and $U_iK$ have the same
capacity and are centrally symmetric, we see that if $l_i$ is a
capacity carrier of $U_iK$ and is homothetic to a capacity carrier
$l_j$ of $U_jK$, then they must actually be identical. We thus
conclude that equality in Corollary~\ref{Corr-Mean-Width} implies in
fact that for every $N$ and $U_1, \ldots, U_N \in U(n)$ we have that
\begin{equation} \label{eq1-eq-in-Mstar} U_1 L \cap \cdots \cap U_NL \neq \emptyset.
\end{equation}

For any $K$ satisfying that $K \neq RB_2^{2n}$ for any $R>0$  there
exists some $\delta>0$ such that for every $\delta$-net ${\cal N}$
on $\sph$ we have that  the restriction of $\| \cdot \|_K$ on ${\cal
N}$ is not constant.  
Assume by contradiction that $K$ satisfies the equality in
Corollary~\ref{Corr-Mean-Width} but is not a Euclidean ball. Fix
$\delta$ as above, 
and fix $C$ such that $C^{-1}|x| \le \|x\|_K \le C|x|$ for all $x$.
Take $U_1 , \ldots, U_N$ to be a $\delta$-net on $U(n)$ with respect
to, say, the operator norm. The finiteness of $N$ follows from
compactness on $U(n)$. Thus for every $U \in U(n)$ there is some $j$
such that $|U_jx - Ux|\le \delta |x|$ for all $x$.

It follows from~$(\ref{eq1-eq-in-Mstar})$ that there exists $l \in
U_1^{-1} L \cap \cdots \cap U_N^{-1}L$. In particular, $U_i l
\subset
\partial K$, and so $\{U_i l (0)\}_{i=1}^N \subset \partial K$.
Consider the set ${\cal N} = \{\frac{U_i l (0)}{|U_i l
(0)|}\}_{i=1}^N \subset \sph$. Note that ${\cal N}$ is a
$\delta$-net of $\sph$.
However, on this set the norm is constant and equals $1/|l(0)|$ (
since $U_i (l(0)) \in \partial K$ so $\| U_i(l(0)) \|_K = 1$).
This is a contradiction to the choice of $\delta$, and we conclude
that the norm $\|\cdot\|_K$ must have been Euclidean, completing the
proof of Corollary~\ref{Corr-Mean-Width}. $\hfill \square$

\vskip 12pt

\noindent {\bf Proof of  Corollary~\ref{Corr-intersections}.} Let
$K,T \subset {\cal K}^{2n}$ be general convex bodies, $x,y \in
\R^{2n}$ and $0\le \lambda \le 1$. It is easy to verify that
\[ \lambda (K \cap (x+T)) + (1-\lambda) (K \cap
(y+T)) \subset (K \cap (\lambda x+(1-\lambda) y+ T)).\]
Therefore, by monotonicity, we have that
\[ c^{1/2}\left(\lambda (K \cap (x+T)) + (1-\lambda) (K \cap
(y+T))\right) \le  c^{1/2}\left(K \cap (\lambda x+(1-\lambda) y+
T)\right).\] The Brunn-Minkowski inequality then implies that
\[
c^{1/2}(\lambda(K \cap (x+T))) +   c^{1/2}((1-\lambda)(K \cap
(y+T))) \le c^{1/2}(K \cap (\lambda x+(1-\lambda) y+ T)),\] and
using homogeneity of capacity the proof of the general case is
complete:
\[ \lambda c^{1/2}(K \cap (x+T)) + (1-\lambda)  c^{1/2}(K \cap
(y+T)) \le c^{1/2}(K \cap (\lambda x+(1-\lambda) y+ T)).\] For the
symmetric case, let $y = -x$ and $\lambda = 1/2$, we get
\[ (1/2) c^{1/2}(K \cap (x+T)) + (1/2)  c^{1/2}(K \cap
(-x+T)) \le c^{1/2}(K \cap T).\] The second term on the left hand
side equals to $c^{1/2}(-K \cap (x-T))$ (since $-Id$ is a symplectic
map), which, by the symmetry assumptions on $K$ and $T$, is the same
as $c^{1/2}(K \cap (x+T))$, the first term, and the inequality
\[  c^{1/2}(K \cap (x+T))  \le c^{1/2}(K \cap T)\]
is established, $\hfill \square$

\section{Proof of Lemma~\ref{tech_lemma}}\label{prooflemma}

The proof is divided into three steps. We follow closely the
arguments in~\cite{HZ} and~\cite{MZ}.

\noindent{\bf First step:} The functional $I$  is bounded from below
on ${\cal E}$. Indeed, the function $h_K^p$ being continues and
homogeneous of degree $p
> 1$ satisfies $$ {\frac 1 \alpha} |y|^{p} \leq h_K^p(y) \leq \alpha
|y|^{p},$$ for some constant $ \alpha \geq1$, and thus
$$ I(z) =  \int_0^{2\pi} h_K^p ( \dot{z}(t)) dt \geq {\frac 1 \alpha} \| \dot{z} \|^{p}_{p},$$
where $\| \cdot \|_p$ stands for the $L_p$ norm on $S^1$. From
H\"{o}lder's inequality it follows that for $z \in {\cal E}$
\begin{equation} \label{eq8} 2 = \int_0^{2\pi} \langle Jz(t),
\dot{z}(t) \rangle dt \leq  \| z \|_{q} \|\dot{z} \|_{p}, \ \ {\rm
where} \ \ {\frac 1 p} + {\frac 1 q} = 1 \end{equation} Using
Poincar\'e inequality and the fact that $ \int_0^{2 \pi} z(t)dt =
0$, we deduce that
\begin{equation} \label{eq9}  \|z \|_{q} \leq \|z \|_{\infty}
\leq 2 \pi \| \dot{z} \|_1 \leq 2 \pi \| \dot{z} \|_{p},
\end{equation} and hence ${
\frac 1 {\sqrt{\pi}}} \leq \| \dot{z} \|_{p}$, which in turn implies
that
\begin{equation}  \label{eq10} I(z) = \int_0^{2\pi} h_K^p ( \dot{z}(t)) dt \geq {\frac 1 \alpha} \|
\dot{z} \|^{p}_{p} \geq  {\frac {{\pi}^{\frac {-p} 2}} \alpha }  > 0
\end{equation}

\noindent{\bf Second step:} The functional $I$ attains its minimum
on ${\cal E}$ namely, there exists $\tilde z \in {\cal E}$ with
$$ \int_0^{2\pi} h_K^p ( \dot{\tilde z}(t)) dt = \inf_{z \in
{\cal E}} \int_0^{2\pi} h_K^p ( \dot{z}(t)) dt = \tilde \lambda >
0$$ To show this, we pick a minimizing sequence $z_j \in {\cal E}$
such that
$$ \lim_{j \rightarrow \infty} \int_0^{2\pi} h_K^p ( \dot{z_j}(t))dt = \tilde \lambda$$
It follows from~$(\ref{eq8})$,~$(\ref{eq9})$, and~$(\ref{eq10})$
that there exists some constant $C > 0$ such that $$ {\frac 1 C}
\leq \| \dot{z}_j \|_{p} \leq C$$ Moreover, from~$(\ref{eq9})$ it
follows that
\begin{equation} \label{eq30} \| z_j \|_p \leq  \|  \dot{z}_j \|_{p} \leq C \end{equation} In
particular, $ z_j $ is a bounded sequence in the Banach space
$W^{1,p}(S^1,{\mathbb R}^{2n})$ and therefore, a subsequence, also
denoted by $ z_j $ converges weakly in $W^{1,p}(S^1,{\mathbb
R}^{2n})$ to an element $z_* \in W^{1,p}(S^1,{\mathbb R}^{2n})$.
Indeed, the closed unit ball of a reflexive Banach space is weakly
compact and the space $W^{1,p}(S^1,{\mathbb R}^{2n})$, where $p >1$,
is known to be reflexive (see e.g.,~\cite{Ad}). We will show below
that $z_* \in {\cal E}$. First we claim that $z_j$ converges
uniformly to $z_*$ i.e.
\begin{equation} \label{conv_unif}
 \sup_t | z_j(t) - z_*(t) | \rightarrow 0
\end{equation}
Indeed, the $z_j$ are uniformly continuous:
$$ | z_j(t) - z_j(s)| \leq  | \int_s^t {\dot z}_j(\tau) d\tau| \leq
|t-s|^{1/q} C ,$$ and the claim follows from the Arzel\`a-Ascoli
theorem. Next we claim that $z_* \in {\cal E}$. Indeed, even the
weak convergence immediately implies that the mean value of $z_*$
vanishes. Moreover,
$$ 2 = \int_0^{2 \pi} \langle J z_j(t) ,  \dot{z}_j(t) \rangle dt =
 \int_0^{2 \pi} \langle J (z_j(t)-z_*(t)) ,  \dot{z}_j(t) \rangle dt +  \int_0^{2 \pi} \langle J z_*(t) ,
 \dot{z}_j(t) \rangle dt$$ The first term on the right hand side tends to zero
 by equation~$(\ref{conv_unif})$, H\"{o}lder inequality,  and
 equation~$(\ref{eq30})$. The second term converges because of the
 weak convergence to $$ \int_0^{2 \pi} \langle J z_*(t) ,
 \dot{z}_*(t)
 \rangle dt$$
To see this, one must check that the linear functional $f(w) =
\int_0^{2 \pi} \langle J z_*(t), \dot w(t) \rangle dt$ is bounded on
$W^{1,p}(S^1, \mathbb R^{2n})$. This follows from H\"{o}lder's
inequality as $z_* \in L_q(S^1, {\mathbb R}^{2n})$ (since $z_* \in
W^{1,p}(S^1, \mathbb R^{2n})$). Thus the equation above, taking
limit $j \rightarrow \infty$ takes the form
$$ \int_0^{2 \pi} \langle J z_*(t) , \dot{z}_*(t)
 \rangle  dt = 2,$$
which implies that $z_* \in {\cal E}$. We now turn to show that $z_*
\in {\cal E}$ is indeed the required minimum. We use the convexity
of $h_K^p$ and deduce the point-wise estimate
$$ \langle \nabla h_K^p({\dot z}_j(t)), {\dot z}_*(t)-{\dot z}_j(t) \rangle \leq h_K^p({\dot z}_*(t)) -
h_K^p({\dot z}_j(t)) \leq \langle \nabla h_K^p({\dot z}_*(t)), {\dot
z}_*(t)-{\dot z}_j(t) \rangle,$$ which gives
\begin{equation} \label{eq1000} \int_0^{ 2 \pi} h_K^p ({\dot z}_*(t))dt - \int_0^{ 2
\pi} h_K^p ({\dot z}_j(t))dt \leq \int_0^{ 2 \pi} \langle \nabla
h_K^p({\dot z}_*(t)), {\dot z}_*(t)-{\dot z}_j(t) \rangle dt
\end{equation}
To see that the right hand side of inequality~$(\ref{eq1000})$ tends
to zero, it is enough as before to check that
 $\nabla h_K^p({\dot z}_*)$ belongs to $L_q(S^1,{\mathbb
R}^{2n})$. Indeed, since $\nabla h_K^p$ is homogeneous of degree
$p-1$ there exists some positive constant $K$ for which $|\nabla
h_K^p(x)| \leq K |x|^{p-1}$, and hence it follows from
equation~$(\ref{eq30})$ that
$$ \int_0^{2 \pi} |\nabla h_K^p({\dot z}_*(t))|^qdt =  \int_0^{2 \pi}
|\nabla h_K^p({\dot z}_*(t))|^{\frac {p} {p-1}}dt \leq K^{\frac {p}
{p-1}} \int_0^{2 \pi} |{\dot z}_*(t) |^p dt < \infty.$$ Thus the
right hand side of inequality~$(\ref{eq1000})$ tends to zero. Hence,
$$ \tilde \lambda \leq \int_{0}^{2 \pi} h_K^p( {\dot z}_*(t)) dt
\leq \liminf_{j \rightarrow \infty} \int_{0}^{2 \pi} h_K^p( {\dot
z}_j(t)) dt = \tilde \lambda,$$ and we have proved that $z_*$ is the
minimum of $I(z)$ for $z \in {\cal E}$.

\noindent{\bf Third step:} First we show that the critical points of
$I$ satisfy the required Euler equation~$(\ref{Euler_eq})$. Let $z$
be a critical point of $I$. Hence, according to
Definition~\ref{def-crit-pt},  for every $\xi \in
W^{1,p}(S^1,\mathbb R^{2n})$ such that $\int_0^{2 \pi} \xi(t)dt =
0$,
$ \int_0^{2 \pi} \langle J z(t) , {\dot \xi}(t) \rangle dt =0$
we have that $$ \int_0^{2 \pi} \langle \nabla h_K^p({\dot z}(t)) ,
{\dot \xi}(t) \rangle dt = 0. $$ Next we choose a special $\xi$
namely such that $\dot \xi$ is of the form $ {\dot \xi} = \nabla
h_K^p({\dot z})- \beta J z - \alpha$ where $\alpha$ is a vector and
$\beta$ is a constant. The vector $ \alpha$ is chosen so that
$\xi(0) = \xi(2 \pi)$ namely $$ \alpha = {\frac 1 {2 \pi}} \int_0^{2
\pi} \nabla h_K^p( \dot z(t))dt. $$ In order to show that $$ x(t) =
\int_0^t \dot \xi(s)ds \in W^{1,p}(S^1,R^{2n}),$$ one uses a simple
continuity properties of $\dot \xi$. We choose $\beta$ so that the
condition $$ \int_0^{2 \pi} \langle Jz(t), \dot \xi(t) \rangle dt=
0$$ is satisfied. With this function $\xi$ we compute
$$ \int_0^{2 \pi} | {\dot \xi}(t) |^2 dt =  \int_0^{2 \pi} \langle \nabla
h_K^p({\dot z}(t)) , {\dot \xi}(t) \rangle dt - \beta \int_0^{2
\pi} \langle Jz(t) , {\dot \xi}(t) \rangle dt - \langle \alpha,
\int_0^{2 \pi} {\dot \xi}(t)dt \rangle  = 0
$$
Thus, the critical point $z$ satisfies the Euler equation $\nabla
h_K^p({\dot z}) =  \beta J z + \alpha$. Moreover, it follows from
Euler formula that
$$ \lambda = \int_0^{2 \pi} h_K^p({\dot z}(t))dt = {\frac 1 p} \int _0^{2 \pi}
\langle  \nabla  h_K^p({\dot z}(t)) , {\dot z}(t) \rangle dt =
{\frac \beta p} \int_0^{2\pi} \langle Jz(t), {\dot z}(t) \rangle
dt = {\frac {2 \beta} p},$$ and hence $\beta = {\frac {\lambda p}
2}$.

For the other direction, namely that any loop $z$ satisfying Euler
equation~$(\ref{Euler_eq})$ is a critical point of $I$, we simply
check that for $\xi \in W^{1,p}(S^1,R^{2n})$ with $\int_0^{2 \pi}
\xi(t)dt =0$ and $\int_0^{2 \pi} \langle Jz(t) , \dot \xi(t)
\rangle dt = 0$ we have
$$\int_0^{2 \pi} \langle \nabla h_K^p (\dot z(t)), \dot \xi(t) \rangle  dt=
\int_0^{2 \pi} \langle {\frac {\lambda p} 2} Jz(t) + \alpha, \dot
\xi(t) \rangle  dt = 0,$$ as required. This concludes the proof of
the Lemma.

\bigskip
\noindent Shiri Artstein-Avidan\\
School of Mathematical Sciences, Tel Aviv University, Tel Aviv 69978, Israel\\
{\it e-mail}: shiri@post.tau.ac.il

\bigskip
\noindent Yaron Ostrover\\
Department of Mathematics, M.I.T, Cambridge MA 02139, USA\\
{\it e-mail}: ostrover@math.mit.edu


\begin{thebibliography}{}

\bibitem{Ad} Adams, R.A. {\it Sobolev Spaces}, Pure and Applied Mathematics, Vol. 65. Academic Press, New York-London, 1975.

\bibitem{AO} Artstein-Avidan, S., Ostrover Y. {\it On Symplectic
Capacities and Volume Radius}, Preprint math.SG/0603411

\bibitem{AMO}Artstein-Avidan, S., Milman, V., Ostrover, Y. {\it The M-ellipsoid,
Symplectic Capacities and Volume}, Preprint math.SG/0604434

\bibitem{AE} Aubin, J.P., Ekeland, I. {\it Applied nonlinear analysis}. Pure and Applied Mathematics, New York, 1984

\bibitem{B} Borell, C. {\it Capacitary inequalities of the
Brunn-Minkwoski type}, Math. Ann. 263 (1983), 179-184.

\bibitem{BL} Brascamp, H.J., Lieb,  E.H. {\it On extensions of the
Brunn-Minkowski and Pr\'{e}kopa-Leindler theorem, including
inequalities for log-concave functions, and with an application to
the diffusion equation}  J. Funct. Anal., 22  (1976), 366-389.

\bibitem{CHLS} Cieliebak, K., Hofer, H., Latschev, J., Schlenk
F. {\it Quantitative symplectic geometry.} math.SG/0506191.

\bibitem{C} Clarke, F. {\it A classical variational principle for
periodic Hamiltonian trajectories}, Proc. Amer. Math. Soc.,
76:186-188, 1979.

\bibitem{CE} Clarke, F., Ekeland, I. {\it Hamiltonian trajectories
having prescribed minimal period}, Comm. Pure Appl. Math. 33, 103-116 (1980)

\bibitem{deG} de Gosson, Maurice. {\it Symplectic geometry and quantum mechanics}, Operator Theory: Advances and Applications, 166.
Advances in Partial Differential Equations (Basel). Birkh\"auser Verlag, Basel, 2006.


\bibitem{Ek} Ekeland, I. {\it Convexity Methods in Hamiltonian
Mechanics}, Ergebnisse der Mathematik und ihrer Grenzgebiete (3), 19. Springer-Verlag, Berlin, 1990.

\bibitem{EH} Ekeland, I., Hofer, H. {\it Symplectic topology
and Hamiltonian dynamics}, Math. Z. 200 (1989), no. 3,
355--378.

\bibitem{EH1} Ekeland, I., Hofer, H. {\it Symplectic topology
and Hamiltonian dynamics II}, Math. Z.  203 (1990), no.4, 553--567.

\bibitem{F} Firey, Wm.J. {\it p-Means of convex bodies}, Math. Scand. 10 (1962), 17-24.

\bibitem{Gar} Gardner, R.J. {\it The Brunn-Minkowski inequality}, Bull. Amer. Math. Soc. 39 (2002), 355-405.

\bibitem{G} Gromov, M. {\it Pseudoholomorphic curves in symplectic manifolds},
Invent. Math. 82 (1985), no. 2, 307-347.

\bibitem{Ho} Hofer, H. {\it Symplectic capacities.} Geometry of
low-dimensional manifolds, 2 (Durham, 1989), 15-34, London Math.
Soc. Lecture Note Ser., 151, Cambridge Univ. Press, Cambridge,
1990.

\bibitem{HZ} Hofer, H., Zehnder, E. {\it Symplectic Invariants
and Hamiltonian Dynamics}, Birkh\"auser, Basel (1994).

\bibitem{HZ1} Hofer, H., Zehnder, E. {\it A new capacity for
symplectic manifolds}, Analysis et cetera. Academic press, 1990.
Pages 405-428.

\bibitem{L} Lalonde, F. {\it Energy and capacities in symplectic
topology.} Geometric topology (Athens, GA, 1993), 328-374, AMS/IP
Stud. Adv. Math., {\bf 2.1}, Amer. Math. Soc., Providence, RI, 1997.

\bibitem{MZ} Moser, J., Zehnder, E. {\it Notes on Dynamical Systems},
Courant Lecture Notes in Mathematics, 12. New York University
(2005).

\bibitem{Mc} McDuff, D. {\it Symplectic topology and capacities},
Prospects in mathematics (Princeton, NJ, 1996), 69-81, Amer. Math.
Soc., Providence, RI, 1999.

\bibitem{McS} McDuff, D., Salamon, D. {\it Introduction to
Symplectic Topology}, 2nd edition, Oxford University Press, Oxford,
England (1998).

\bibitem{Pisbook}  Pisier, G. {\it The volume of convex bodies and Banach space geometry},
Cambridge Tracts in Mathematics, 94. Cambridge University Press,
Cambridge, 1989.

\bibitem{R} Rabinowitz, P. {\it Periodic solutions of Hamiltonian
systems}, Comm. Pure Appl. Math., 31:157-184, 1978.

\bibitem{Sch} Schneider, R. {\it Convex bodies: the Brunn-Minkowski theory},
Encyclopedia of Mathematics and its Applications, 44. Cambridge
University Press, 1993.

\bibitem{Sib} Siburg, K.F. {\it Symplectic capacities in two
dimensions}, manusripta math., 78:149-163, 1993.


\bibitem{V} Viterbo, C.
{\it Metric and isoperimetric problems in symplectic geometry.} J.
Amer. Math. Soc. 13 (2000), no. 2, 411--431.

\bibitem{V1} Viterbo, C. {\it Capacit\'es symplectiques
et applications (d'apr\`es Ekeland-Hofer, Gromov).} S\'eminaire
Bourbaki, Vol. 1988/89. Ast\'erisque no.  177-178 (1989), Exp. no.
714, 345-362.

\bibitem{W} Weinstein, A. {\it Periodic orbits for convex
Hamiltonian systems}, Ann. Math., 108:507-518, 1978.

\end{thebibliography}
\end{document}